\documentclass[11pt]{article}
\usepackage[margin=1.1in]{geometry}
\usepackage{amsmath,amssymb,amsthm,mathtools,enumitem}
\usepackage[hidelinks]{hyperref}

\newtheorem{theorem}{Theorem}[section]
\newtheorem{lemma}{Lemma}[section]

\newtheorem{corollary}{Corollary}[section]
\theoremstyle{definition}
\newtheorem{definition}{Definition}[section]
\theoremstyle{remark}

\numberwithin{equation}{section}

\newcommand{\C}{\mathbb{C}}
\newcommand{\Z}{\mathbb{Z}}
\newcommand{\Q}{\mathbb{Q}}
\newcommand{\RR}{\mathbb{R}}
\newcommand{\N}{\mathcal{N}}
\newcommand{\tr}{\operatorname{tr}_G}
\newcommand{\dimG}{\dim_G}
\newcommand{\Fix}{\operatorname{Fix}}
\newcommand{\Rose}{R} 

\newcommand{\ip}[2]{\langle #1,#2\rangle}
\newcommand{\ran}{\operatorname{ran}}
\newcommand{\Vqf}{V_{qf}}

\title{Chain flaring fails on the rose for atoroidal irreducible automorphisms}
\author{Jorge Andres Robinson Arrieta\\[2pt]
\small\texttt{jorge.robinson@keiseruniversity.edu}}
\date{}

\begin{document}
\maketitle

\begin{abstract}
We show that the  chain flare condition introduced by  Clay fails for every topological representative on the topological $n$-rose of every power of an atoroidal irreducible automorphism of $F_n$.

\end{abstract}

\medskip


\section{Introduction}\label{sec:intro}

Let $F=F_n$ be the free group  $\langle x_1,\dots,x_n\rangle$ with $n\ge 2$, and let $\Phi\in\operatorname{Aut}(F)$ be an automorphism of $F$.  The free-by-cyclic group associated with $\Phi$ is the group $G_{\Phi} = F\rtimes_{\Phi} \Z$ with presentation 

\begin{equation}\label{eq:G}
G_{\Phi}=\langle F,t\mid t^{-1}xt=\Phi(x)\text{ for }x\in F\rangle.
\end{equation}
Write $\chi\colon G_{\Phi}\to\Z$ for the canonical epimorphism given by
$\chi(x)=0$ for $x\in F$ and $\chi(t)=1$.

 The $L^2$-torsion, denoted $\rho^{(2)}(G)$, is an analytic invariant introduced in analytic form by Lott \cite{Lott1992} and Mathai \cite{Mathai1992} and in combinatorial form by Carey–Mathai \cite{CareyMathai1992} and Lück–Rothenberg \cite{LueckRothenberg1991} for suitable $L^2$-acyclic groups, including free-by-cyclic groups \cite[Chapters 3 and 13]{Lu02}. Clay showed in \cite{Cl17} that $\rho^{(2)}(G_{\Phi})$ vanishes whenever $\Phi$ is polynomially growing. L\"uck--Schick proved in \cite{LS99} that if $M$ is a compact connected $3$-manifold whose interior carries a complete hyperbolic metric of finite volume, then 
\[
\rho^{(2)}(\pi_1(M))=-\frac{\operatorname{vol}(M)}{6\pi}.
\]
Thus $-\rho^{(2)}(\pi_1(M))>0$. Motivated by the $3$-manifold case, Clay conjectured \cite[Section 1]{Cl17} that if $G$ is a free-by-cyclic group with exponentially growing monodromy, then $-\rho^{(2)}(G)>0$. Clay introduced the chain flare condition in  order to study the positivity of $-\rho^{(2)}(G)$ for free-by-cyclic groups with exponentially growing monodromy. He proved the following result, which we state in the case where $\Phi$ is fully irreducible.

\begin{theorem}[Clay {\cite[Theorem~1.1]{Cl21}}]\label{thm:clay}
Let $\Phi\in\operatorname{Aut}(F_n)$ be fully irreducible. If $f_{\Phi}$ satisfies the chain flare condition relative to the empty subgraph, then
\[
-\rho^{(2)}(G_\Phi)=\int_{1<|z|}\log|z|\,d\mu_L>0 ,
\] where $\mu_L$ is the Brown measure of $L$ introduced in \cite{Br86}.
\end{theorem}


An automorphism $\Phi$ is atoroidal if it has no nontrivial periodic conjugacy class. An automorphism $\Phi$ is fully irreducible  if it has no periodic conjugacy classes of proper nontrivial free factors in $F$. Dowdall–Kapovich–Leininger proved \cite[Corollary~B.4]{DKL2015} that if $\Phi$ is atoroidal and irreducible, then it is fully irreducible.
In this paper, we prove that the chain flare condition fails  for every topological representative on the topological $n$-rose of every power of an atoroidal irreducible automorphism of $F_n$

We now introduce some notation in order to state the main results. For $X\in M_{n\times m}(\C[G])$, $R_X\colon\ell^2(G)^n\to\ell^2(G)^m$ is the $G$-equivariant operator given by $(R_X\xi)_j=\sum_i\xi_i\cdot X_{ij}$. Consider the operator $L:=R_{tJ}$ where $J_{ij}=\partial\Phi(x_i)/\partial x_j$ is the matrix of Fox derivatives of $\Phi$. In \cite[Theorem 4.10]{Cl17} Clay showed that $-\rho^{(2)}(G_{\Phi})=\log\det_{G}(I-L)$. 
We say that $L$ has a spectral gap if there is $r>1$ such that $\sigma(L)\cap\{w\in\C: r^{-1}<|w|<r\}=\emptyset$.


\begin{theorem}\label{cor:B}
$I-L$ is not invertible, and in fact not surjective, on $\ell^2(G)^n$.
\end{theorem}



\begin{theorem}\label{thm:D}
Let $H\subset\Rose_n$ be the empty subgraph or the vertex $*$. The following are equivalent.
\begin{enumerate}[label=\textup{(\alph*)}]
\item $L$ has a spectral gap.
\item $\sigma(L)\cap S^1=\emptyset$.
\item $I-L$ is invertible.
\item For some $N\ge1$, $f_{\Phi^N}$ satisfies Clay's chain flare condition relative to $H$ with $\Vqf=\{0\}$ (Definition~\ref{def:CF}).
\item For some $N\ge1$ and $\lambda>1$, $\lambda\|x\|\le\max\{\|A^Nx\|,\|A^{-N}x\|\}$ for all $x\in\ell^2(F)^n$.
\end{enumerate}
Moreover, if $f_{\Phi}$ satisfies the chain flare condition relative to $H$ with $\Vqf=\{0\}$ and constant $\lambda$, then $\sigma(L)\cap\{\lambda^{-1}<|w|<\lambda\}=\emptyset$.
\end{theorem}

By Theorem \ref{cor:B}, none of (a)--(e) in Theorem \ref{thm:D} ever holds.

\begin{corollary}\label{cor:E}
Let $\Psi\in\operatorname{Aut}(F_n)$ and $N\ge 1$. If $\Fix(\Psi^N)=\{1\}$, then $f_{\Psi^N}$ does not satisfy Clay's chain flare condition relative to  $H\in\{\emptyset,*\}$.

\end{corollary}

\begin{corollary}\label{cor:F}
Let $\Phi\in\operatorname{Aut}(F_n)$ be atoroidal and irreducible. Then for every $N\ge1$, every proper $f_{\Phi^N}$-invariant subgraph is  $\emptyset$ or *. Thus, the rose map $f_{\Phi^N}$ does not satisfy Clay's chain flare condition relative to  $H\in\{\emptyset,*\}$.
\end{corollary}


\section{Background}\label{sec:background}


We first recall the definition of von Neumann dimension. Let $\N(G)$ be the group von Neumann algebra and $M_n(\N(G))$ denote the algebra of bounded operators  $\ell^2(G)^n$ commuting with the left action of $G$. The trace is defined by $\tr(T):=\sum_i\ip{T\delta_e^{(i)}}{\delta_e^{(i)}}$, for $T\in M_n(\mathcal{N}(G))$. A Hilbert $G$-module is a closed $G$-invariant subspace $U\subseteq\ell^2(G)^n$, and we define $\dimG U:=\tr(P_U)$. A morphism between Hilbert $G$-modules is a bounded $G$-equivariant linear map.  We recall the basic properties of the von Neumann dimension.

\begin{theorem}[L\"uck {\cite[Ch.~1]{Lu02}}]\label{thm:dim}

The following properties are satisfied by the von Neumann dimension:

\begin{enumerate}[label=\textup{(\roman*)}]

\item $\dimG$ is additive on orthogonal sums, $\dimG\ell^2(G)^n=n$, and $\dimG U=0$ only if $U=\{0\}$.
\item For $T\in M_n(\N(G))$, $\dimG\ker T=\dimG\ker T^*$.
\item If $T\in M_n(\N(G))$ is bounded below, then $T$ is invertible and $T^{-1}\in M_n(\N(G))$.
\item A bijective morphism $U\to V$ of Hilbert $G$-modules has a bounded equivariant inverse, and $\dimG U=\dimG V$.
\end{enumerate}
\end{theorem}

 For a morphism $T\colon U\to V$ of Hilbert $G$-modules, $F_T(z)=\tr\big(E^{T^*T}_{z^2}\big)$ is the spectral density function, where $E^{T^*T}_z$ is the spectral family of $T^*T$, and the Fuglede--Kadison determinant is
\[
\det\nolimits_G(T)=\exp\int_{0+}^\infty\log z\,dF_T(z)\quad\text{if }\int_{0+}^\infty\log z\,dF_T(z)>-\infty,\qquad \det\nolimits_G(T)=0\ \text{otherwise};
\] see \cite[Def.~3.11]{Lu02}. 



We next recall Fox calculus and the rose-map chain operator.

For each $j$ there is a unique map $\partial/\partial x_j\colon F\to\Z[F]$ with $\partial(gh)/\partial x_j=\partial g/\partial x_j+g\,\partial h/\partial x_j$ and $\partial x_i/\partial x_j=\delta_{ij}$; it extends linearly to $\Z[F]$. We recall the fundamental formula of Fox calculus \cite{Fox53}:

\begin{theorem}[Fox {\cite[(2.3)]{Fox53}}]\label{thm fundamental fox}
For every $x\in F$,
\begin{equation}\label{eq:fox}
  \Phi(x)-1=\sum_j \frac{\partial\Phi(x)}{\partial x_j}\,(x_j-1).
\end{equation}
\end{theorem}

Suppose $f : \Gamma \xrightarrow{} \Gamma$ is a homotopy equivalence that fixes a vertex $*\in V(\Gamma)$ and that $H \subset \Gamma$ is an $f$-invariant subgraph. Fixing an isomorphism $\pi_1(\Gamma, *) \cong F$ we have that $f$ and the trivial path based at $*$ induce an automorphism of $F$ that we denote by $\Phi_f$ .
We can lift the topological representative
$f:\Gamma\to\Gamma$ of $\phi=[\Phi_f]\in\operatorname{Out}(F)$
to a cellular map
$\widetilde{f}:\widetilde{\Gamma}\to\widetilde{\Gamma}$,
where $\widetilde\Gamma$ is the universal cover of $\Gamma$. Fix a lift $\tilde{*}$ of $*$ and let $\tilde f$ be the lift with $\tilde f(\tilde{*}) = \tilde{*}$.
Although $\widetilde f$ is not $F$-equivariant, it satisfies
\[
\widetilde f(gz)
=
\Phi_f(g)\widetilde{f}(z),
\]
for every $g\in F$ and $z\in\widetilde\Gamma$,
where $\Phi_f\in\operatorname{Aut}(F)$ represents $\phi$.
The induced map on cellular $1$-chains is denoted by
\[
A_f:C_1(\widetilde\Gamma;\Q)
\longrightarrow
C_1(\widetilde\Gamma; \Q).
\] The inverse of $A_f$ is given by $A_h$, where $h:\Gamma\xrightarrow{}\Gamma$ is a topological representative of $\Phi_f^{-1}$.




We recall the Riesz decomposition theorem.

\begin{theorem}
    [{Riesz decomposition; \cite[Thm.~III.6.17]{Ka95}}]\label{thm:riesz}
Let $T$ be a bounded operator on a Banach space $X$ such that $\sigma(T)\cap S^1=\emptyset$, and write $\sigma(T)=\sigma_{s}\cup\sigma_{u}$ with $\sigma_{s}=\sigma(T)\cap\{|w|<1\}$  and $\sigma_{u}=\sigma(T)\cap\{|w|>1\}$. Let 
\[
P=\frac1{2\pi i}\oint_{S^1}(w-T)^{-1}\,dw .
\]
Then $P$ is a bounded idempotent commuting with $T$; the subspaces $M_{s}=\ran P$ and $M_{u}=\ran(I-P)$ are closed and $T$-invariant, $X=M_{s}\oplus M_{u}$ (a direct sum, not necessarily orthogonal); and $\sigma(T|_{M_{s}})=\sigma_{s}$, $\sigma(T|_{M_{u}})=\sigma_{u}$.
\end{theorem}{}

We also recall Clay's chain flare condition  for a topological representative $f$ in the case $H\in \{\emptyset,*\}$, and  $\Vqf=\{0\}$ \cite[(CFH1)-- (CFH3)]{Cl21}. Let $C_1(\tilde{\Gamma}, \tilde{H}; \mathbb{Q})$ denote the module of finitely supported rational 1-chains, with the norm of $\ell^2(F)^n$.

\begin{definition}[{\cite[\S2.4]{Cl21}}]\label{def:CF}
The map $f$ satisfies the \emph{chain flare condition relative to $H\in\{\emptyset,*\}$ with $\Vqf=\{0\}$} and constant $\lambda>1$ if 
\[
\lambda\|x\|\le\max\{\|A_fx\|,\|A_f^{-1}x\|\}\qquad\text{for all }x\in C_1(\tilde{\Gamma}, \tilde{H}; \mathbb{Q}) .
\]
\end{definition}

\section{The intertwining relation and spectral properties of $R_t$}\label{sec:AB}

Let $d=(x_1-1,\dots,x_n-1)^{T}$. We have the following basic lemmas.

\begin{lemma}For composable $X$ and $Y$
\begin{equation}\label{eq:comp}
R_Y\circ R_X=R_{XY}.
\end{equation}
\begin{proof}
    
$(R_Y(R_X\xi))_k=\sum_j(\sum_i\xi_i\cdot X_{ij})\cdot Y_{jk}=\sum_i\xi_i\cdot(\sum_jX_{ij}Y_{jk})=(R_{XY}\xi)_k$.

\end{proof}
\end{lemma}

\begin{lemma}\label{lem:intertwine}
$(tJ)\,d=d\,t$ in $M_{n\times1}(\Z[G])$. Consequently $R_d\circ L=R_t\circ R_d$.
\end{lemma}
\begin{proof}
The $i$-th entry of $(tJ)d$ is $\sum_jtJ_{ij}(x_j-1)=\sum_jt\frac{\partial \Phi(x_i)}{\partial x_j}(x_j-1)$.  By Theorem $\ref{thm fundamental fox}$ this is equal to  $t(\Phi(x_i)-1)=(x_i-1)t$ since $t\Phi(x_i) = x_it$ by (1.1). Thus, by (3.1), $R_d\circ L=R_d\circ R_{tJ}=R_{(tJ)d}=R_{dt}=R_t\circ R_d$.
\end{proof}

\begin{lemma}\label{lem:surj}
$R_d\colon\ell^2(G)^n\to\ell^2(G)$ is surjective.
\end{lemma}
\begin{proof}
Notice that $d^*d=\sum_k(x_k^{-1}-1)(x_k-1)=\sum_k(2-x_k-x_k^{-1})=:\Delta_0$.  Hence $R_d\circ R_{d^*}=R_{d^*d}=R_{\Delta_0}$ is invertible by \cite[Theorem 3]{Ke59}.
\end{proof}

\begin{lemma}\label{lem:Rt} $R_t$ satisfies the following properties:
\begin{enumerate}[label=\textup{(\roman*)}]

\item $R_t\in\N(G)$ is unitary.
\item $R_t-z$ is injective for every $z\in\C$.
\item $\sigma(R_t)=S^1$, and $I-R_t$ is not bounded below.
\end{enumerate}
\end{lemma}
\begin{proof}
(i) $R_t^*=R_{t^{-1}}=R_t^{-1}$. (ii) For $z=0$ use (i). Let $z\ne0$ and suppose that $\xi\in\ker(R_t-z)$. Then $\xi t=z\xi$ and $\xi_{gt^m}=z^{-m}\xi_g$ for all $g\in G$, $m\in\Z$. Since $|\xi_g|^2\sum_{m\in\Z}|z|^{-2m}\le\|\xi\|^2$ and $\sum_{m\in\Z}|z|^{-2m}$ diverges, $\xi_g=0$. (iii) For $|z|=1$ the unit vectors $\xi_N=N^{-1/2}\sum_{k=0}^{N-1}z^{-k}\delta_{t^k}$ satisfy $\|(R_t-z)\xi_N\|=(2/N)^{1/2}$, so $S^1\subseteq\sigma(R_t)$ and $I-R_t$ is not bounded below. The containment $\sigma(R_t)\subseteq S^1$ follows since $R_t$ is unitary. 



\end{proof}


\begin{proof}[Proof of Theorem~\ref{cor:B}]
By Lemma~\ref{lem:intertwine}, $(I-tJ)d=d(1-t)$, so $R_d(I-L)=(I-R_t)R_d$. Suppose $I-L$ is surjective. Given $\zeta\in\ell^2(G)$, choose $\eta$ with $R_d\eta=\zeta$ (Lemma~\ref{lem:surj}) and $\xi$ with $(I-L)\xi=\eta$; then $\zeta=(I-R_t)R_d\xi$. So $I-R_t$ is surjective; it is injective by Lemma~\ref{lem:Rt}(ii), hence invertible by the bounded inverse theorem, hence bounded below, contradicting Lemma~\ref{lem:Rt}(iii).
\end{proof}

\section{Spectral gap and chain flaring: proof of Theorem~\ref{thm:D}}\label{sec:D}

\begin{lemma}\label{lem:shift}
We have an orthogonal direct sum decomposition $\ell^2(G)^n=\bigoplus_{\ell\in\Z}t^\ell\ell^2(F)^n$, and $L(\xi)=\sum_{\ell\in\Z}t^{\ell+1} A(\xi^{(\ell)})$ for $\xi=\sum_{\ell\in\Z}t^{\ell}\xi^{(\ell)}\in\ell^2(G)^n$ with $\xi^{(\ell)}\in \ell^{2}(F)^{n}$. Hence $L$ is invertible and $\|L^{N}\eta\|=\|A^{N}\eta\|$ for $\eta\in\ell^2(F)^n$, $N\in\Z$.
\end{lemma}
\begin{proof}
The fact that $G=\bigsqcup_\ell t^{\ell}F$ gives the decomposition. The formula (7.1) in \cite{Cl21} gives the formula for $L$. The inverse of $L$ is given by  $L^{-1}(\xi)=\sum_{\ell\in\Z}t^{\ell-1} A^{-1}(\xi^{(\ell)})$.
\end{proof}




\begin{lemma}\label{lem:cxflare}
Let $T$ be invertible on a Banach space and $\kappa>1$ with $\kappa\|x\|\le\max\{\|Tx\|,\|T^{-1}x\|\}$ for all $x$. If $T=S^k$ with $k\ge1$, then $S-v$ is bounded below for $\kappa^{-1/k}<|v|<\kappa^{1/k}$.
\end{lemma}
\begin{proof}

Let $u=v^k$, so $\kappa^{-1}<|u|<\kappa$, and put $\mu=\max\{|u|,|u|^{-1}\}<\kappa$. For a unit vector $x$, let $e:=\|(T-u)x\|$. Then we have $\|Tx\|\le \|ux\|+\|(T-u)x\| \le|u|+e$. Write $T^{-1}x=u^{-1}x-u^{-1}T^{-1}(T-u)x$. Then $\|T^{-1}x\|\le|u|^{-1}+|u|^{-1}\|T^{-1}\|e$. Then $\kappa\le\max\{\|Tx\|,\|T^{-1}x\|\}\le \max\{|u|+e,|u|^{-1}+|u|^{-1}\|T^{-1}\|e\}\le\mu+ce$ with $c=\max\{1,|u|^{-1}\|T^{-1}\|\}$. Thus, $e=\|(T-u)x\|\ge\dfrac{\kappa-\mu}{c}$, so $T-u$ is bounded below. Since $S^k-v^k=Q(S-v)$ with $Q=\sum_{j<k}v^jS^{k-1-j}$, $S-v$ is also bounded below because $\|(S -v)x\| \ge  \dfrac{(\kappa - \mu)}{c\|Q\|}$ .
\end{proof}

\begin{lemma}\label{lem:realflare}
Let $B$ be real and invertible on $K$, and $\lambda>1$ with $\lambda\|y\|\le\max\{\|By\|,\|B^{-1}y\|\}$ for all $y\in K_\RR$. Then for every $k\ge1$:
\begin{enumerate}[label=\textup{(\roman*)}]
\item $\lambda^k\|y\|\le\max\{\|B^ky\|,\|B^{-k}y\|\}$ for $y\in K_\RR$;
\item $\frac{\lambda^k}{\sqrt2}\|x\|\le\max\{\|B^kx\|,\|B^{-k}x\|\}$ for $x\in K$.
\end{enumerate}
\end{lemma}
\begin{proof}
(i) This is Lemma 2.2(3) in \cite{Cl21}. (ii) Write $x=p+iq$ with $p,q\in K_\RR$; then $\|B^{\pm k}x\|^2=\|B^{\pm k}p\|^2+\|B^{\pm k}q\|^2$. By (i) there are signs $\varepsilon,\delta$ with $\|B^{\varepsilon k}p\|\ge\lambda^k\|p\|$ and $\|B^{\delta k}q\|\ge\lambda^k\|q\|$. If $\varepsilon=\delta$, then $\|B^{\varepsilon k}x\|\ge\lambda^k\|x\|$. Otherwise the larger of $\|B^{\varepsilon k}x\|^2$, $\|B^{\delta k}x\|^2$ is at least $\lambda^{2k}\max\{\|p\|^2,\|q\|^2\}\ge\frac{\lambda^{2k}}2\|x\|^2$.
\end{proof}

\begin{lemma}\label{lem:fourierL}
Let $s,c>0$ with $\|(A-u)y\|\ge c\|y\|$ for all $y\in\ell^2(F)^n$ and $|u|=s$. Then $\|(L-w)\xi\|\ge c\|\xi\|$ for all $\xi\in\ell^2(G)^n$ and $|w|=s$.
\end{lemma}
\begin{proof}
It suffices to treat $\xi=\sum_kt^k\cdot\eta_k$ with finitely many $\eta_k\ne0$. Put $p(z)=\sum_kz^k\eta_k$ for $z\in S^1$. Expanding $\|p(z)\|^2=\sum_{k,l}z^{k-l}\ip{\eta_k}{\eta_l}$ and integrating gives Parseval's identity $\int_{S^1}\|p\|^2dm_{S^1}=\sum_k\|\eta_k\|^2=\|\xi\|^2$. By Lemma~\ref{lem:shift}, $(L-w)\xi=\sum_kt^k\cdot(A\eta_{k-1}-w\eta_k)$, whose transform is $q(z)=zAp(z)-wp(z)=z(A-\bar zw)p(z)$. Since $|\bar zw|=s$, $\|q(z)\|\ge c\|p(z)\|$, and Parseval for $q$ and $p$ gives the claim.
\end{proof}

\begin{proof}[Proof of Theorem~\ref{thm:D}] 

(a)$\Rightarrow$(b)$\Rightarrow$(c) is clear. 

(c)$\Rightarrow$(b). Suppose that $I-L$ is invertible. Let $z$ with $|z|=1$ and let $T_{z}$ be the operator on $\ell^{2}(G)^n$ given by $T_{z}(t^{\ell}x)=z^{\ell}t^{\ell}x$. Then $zI-L=zT_{z^{-1}}(I-L)T_z$, so $zI-L$ is invertible.

(b)$\Rightarrow$(a) since $\sigma(L)$ and $S^1$ are disjoint compact sets.

\emph{(e)$\Rightarrow$(d)}  follows by restricting (e) to rational chains. We now prove (d)$\Rightarrow$(e)$\Rightarrow$(a)$\Rightarrow$(e). 

\emph{(d)$\Rightarrow$(e).} 
This follows from Lemma \ref{lem:realflare}.

\emph{(e)$\Rightarrow$(a).} By Lemma~\ref{lem:cxflare} with $T=A^N$, $S=A$, the operator $A-u$ is bounded below by a positive constant $c(u)=\inf_{\|y\|=1}\|(A-u)y\|$ for $\lambda^{-1/N}<|u|<\lambda^{1/N}$. For such $s=|u|$, continuity of $u\mapsto\inf_{\|y\|=1}\|(A-u)y\|$ and compactness of the circle of radius $s$ give a uniform constant $c_s>0$ with $\|(A-u)y\|\ge c_s\|y\|$. Lemma~\ref{lem:fourierL} gives $\|(L-w)\xi\|\ge c_s\|\xi\|$ for $|w|=s$. By Theorem~\ref{thm:dim}(iii), $L-w$ is invertible. So $\sigma(L)\cap\{\lambda^{-1/N}<|w|<\lambda^{1/N}\}=\emptyset$.


\emph{(a)$\Rightarrow$(e).} Let $\sigma(L)\cap\{r^{-1}<|w|<r\}=\emptyset$ and take the Riesz decomposition $\ell^2(G)^n=M_s\oplus M_u$ for $L$ as in Theorem~\ref{thm:riesz}. Since $L$ is invertible and commutes with $P$, $L^{-1}$ preserves both summands.  Write $L_s:=L|_{M_s}$ and $L_u:=L|_{M_u}$. Then $\sigma(L_s)$ and $\sigma(L_u^{-1})$ lie in $\{|w|\le r^{-1}\}$, so for fixed $1<R<r$ there is $C\ge1$ with $\|(L_s)^k\|,\|(L_u^{-1})^k\|\le CR^{-k}$ for $k\ge0$. Let $p=\max\{\|P\|,\|I-P\|\}$, and $x=x_s+x_u$ with $x_s=Px$ and $x_u=(I-P)x$.  Since $L^Nx_{u}\in M_{u}$ and $x_{u}=(L_{u}^{-1})^NL^Nx_{u}$, $\|x_{u}\|\le CR^{-N}\|L^Nx_{u}\|$. Thus, \ $\|L^Nx_{u}\|\ge C^{-1}R^N\|x_{u}\|$, and $\|L^Nx_{s}\|=\|L_{s}^Nx_{s}\|\le CR^{-N}\|x_{s}\|\le Cp\,R^{-N}\|x\|$. Similarly, $\|L^{-N}x_{s}\|\ge C^{-1}R^N\|x_{s}\|$, and $\|L^{-N}x_{\rm u}\|\le Cp R^{-N}\|x\|$. By the triangle inequality we get
\[
\|L^Nx\|\ge C^{-1}R^N\|x_u\|-Cp\,R^{-N}\|x\|,\qquad\|L^{-N}x\|\ge C^{-1}R^N\|x_s\|-Cp\,R^{-N}\|x\|,
\]
and one of $\|x_s\|,\|x_u\|$ is at least $\|x\|/2$. Hence $\max\{\|L^Nx\|,\|L^{-N}x\|\}\ge\beta_N\|x\|$ with $\beta_N=R^N/(2C)-Cp\,R^{-N}$. Take $\lambda=\beta_N\ge2$ for $N$ big enough. Restricting to $x\in\ell^2(F)^n$ and using Lemma~\ref{lem:shift} gives (e).
\end{proof}

\section{Consequences}\label{sec:consequences}

\begin{proof}[Proof of Corollary~\ref{cor:E}]


Let $\tilde f$ be the lift of $f_{\Psi^N}$ with $\tilde f(\tilde *)=\tilde *$. The vertices fixed by $\tilde f$ correspond to $\Fix(\Psi^N)$. With only one fixed vertex there is no Nielsen $1$-chain, since a Nielsen $1$-chain is
of the form $[u,v]$ with $u\neq v$ vertices fixed by $\tilde f$
\cite[Definition~2.1]{Cl21}. Hence \cite[(CFH3)]{Cl21} forces $V_{qf}=\{0\}$. It follows from Theorem~\ref{thm:D},(d)$\Rightarrow$(c), and Theorem~\ref{cor:B} that the preferred lift of $f_{\Psi^N}$ does not satisfy the Clay's chain flare condition relative to  $H\in\{\emptyset,*\}$.
\end{proof}

\begin{proof}[Proof of Corollary~\ref{cor:F}]

Since $\Phi$ is atoroidal, $\operatorname{Fix}(\Phi^N)=\{1\}$. Since $\Phi$ is fully
irreducible by \cite[Corollary B.4]{DKL2015}, $f_{\Phi^N}$ has no proper invariant subgraph containing an edge: the
petals in such a subgraph would generate a proper nontrivial free factor
$\langle S\rangle$ with $\Phi^N(\langle S\rangle)\subseteq\langle S\rangle$, and a free
factor of $\langle S\rangle$ of full rank is $\langle S\rangle$ itself, so
$\langle S\rangle$ would be $\Phi^N$-invariant. Hence every proper $f_{\Phi^{N}}$-invariant subgraph is $H=\emptyset$ or $H=\{*\}$, and
Corollary~1.1 with $\Psi=\Phi$ applies.
\end{proof}  Every topological representative of $[\Phi^{N}]$ on the $n$-rose is $f_{\Psi}$, where $\Psi \in [\Phi^{N}]$ is induced by the lift fixing $\tilde{*}$, and $\Psi$ is again fully irreducible and atoroidal. Therefore, by Corollary \ref{cor:F} applied to $\Psi$ with $N = 1$, $f_{\Psi}$ does not satisfy Clay's chain flare condition relative to any proper invariant subgraph.

\section*{Declaration on the use of AI}\label{sec:aidecl}
The initial manuscript was developed in conversation with Claude, an AI system made by Anthropic (model Claude Opus~5.5). The author had the intuition that the chain flare condition would be equivalent to $I-L$ having a spectral gap, and asked Claude about this. Claude proposed the main argument for the proof of Theorem \ref{thm:D}. The author directed the investigation, prepared and typed the manuscript, simplified and corrected Claude's arguments, and filled in the gaps. 

\bibliographystyle{abbrv}
\bibliography{chain_flare_references}

\end{document}